\documentclass{amsproc}
\usepackage[T1]{fontenc}
\usepackage{mathtools,amssymb}
\usepackage{mathrsfs}
\usepackage[expansion=false]{microtype}
\usepackage{hyperref}
\usepackage{enumitem}
\hypersetup{hidelinks,pdftitle={A Stinespring-type theorem for complete subnormality preservers},pdfauthor={Pawel Pietrzycki and Jan Stochel}}
\newtheorem{theorem}{Theorem}[section]
\newtheorem{proposition}[theorem]{Proposition}
\newtheorem{lemma}[theorem]{Lemma}

\theoremstyle{definition}

\numberwithin{equation}{section}
\newcommand{\C}{\mathbb C}
\newcommand{\Hh}{\mathcal H}
\newcommand{\Kk}{\mathcal K}
\newcommand{\Ee}{\mathcal E}
\newcommand{\B}{\mathscr B}
\newcommand{\sst}{\sigma\text{-strong}^{*}}
\newcommand{\Phihat}{\widehat\Phi}
\DeclareMathOperator{\ran}{ran}
\setlist[enumerate]{label=\textup{(\roman*)},leftmargin=2.1em,itemsep=2pt,topsep=4pt,parsep=0pt}
\allowdisplaybreaks[1]
\title[A Stinespring-type theorem]{A Stinespring-type theorem for maps\\ completely preserving subnormality}
\author[Pawe\l{} Pietrzycki]{\texorpdfstring{Pawe\l{} Pietrzycki}{Pawel Pietrzycki}}
\address{Institute of Mathematics, Jagiellonian University, \L{}ojasiewicza 6, 30-348 Krak\'ow, Poland}
\email{pawel.pietrzycki@uj.edu.pl}
\author[Jan Stochel]{Jan Stochel}
\address{Institute of Mathematics, Jagiellonian University, \L{}ojasiewicza 6, 30-348 Krak\'ow, Poland}
\email{jan.stochel@uj.edu.pl}
\begin{document}
\begin{abstract}
Stinespring's representation theorem plays a central role in operator
algebras and can be viewed as a solution to the complete preserver
problem for positivity. We replace positivity by subnormality and prove
a corresponding representation theorem for bounded linear maps from
arbitrary $C^*$-algebras into $\B(\Hh)$. Such a map completely preserves
subnormality if and only if it is a simultaneous invariant restriction
of a $*$-representation multiplied by a commuting normal operator.
The normal operator can be chosen as the minimal normal extension of
the canonical weight.
\end{abstract}
\maketitle

\section{Introduction}

For a unital $C^*$-algebra $\mathscr{A}$ and a positive linear functional
$\varphi\colon\mathscr{A}\to\C$, the Gelfand--Naimark--Segal (GNS)
construction~\cite{Segal} gives a $*$-representation
$\pi_\varphi\colon\mathscr{A}\to\B(\Hh_\varphi)$ with cyclic vector
$\xi_\varphi\in\Hh_\varphi$ such that
\[
\varphi(a)=\langle\xi_\varphi,\pi_\varphi(a)\xi_\varphi\rangle
\qquad(a\in \mathscr{A}).
\]
Stinespring's representation theorem extends this construction from
scalar-valued functionals to operator-valued maps, with complete positivity
in place of positivity. Here complete positivity means positivity of every
entrywise matrix amplification.

\begin{theorem}[Stinespring {\cite[Theorem~1]{Stinespring}}]
\label{thm:stinespring}
Let $\mathscr{A}$ be a unital $C^*$-algebra, let $\Hh$ be a Hilbert space, and
let $\Psi\colon \mathscr{A}\to\B(\Hh)$ be a linear map. Then $\Psi$ is completely
positive if and only if there are a Hilbert space $\Kk$, a unital
$*$-representation $\pi\colon \mathscr{A}\to\B(\Kk)$, and a bounded operator
$V\colon\Hh\to\Kk$ such that
\[
\Psi(a)=V^*\pi(a)V\qquad(a\in \mathscr{A}).
\]
In this representation $V^*V=\Psi(1)$. In particular, if $\Psi$ is
unital, then $V$ is an isometry and $\Psi$ is a compression of $\pi$.
\end{theorem}

Arveson's work made this principle a central tool in the study of
nonselfadjoint operator algebras. His extension theorem allows a
completely positive map from an operator system in a unital $C^*$-algebra
to $\B(\Hh)$ to be extended to the whole algebra. Together with
Stinespring's theorem, the corresponding extension result for unital
completely contractive maps on unital subspaces of $C^*$-algebras yields
dilations by $*$-representations
\cite[Theorems~1.2.3 and~1.2.9]{Arveson}.
These tools also underlie his theory of boundary representations:
irreducible representations whose restrictions to a generating unital
subspace have unique completely positive extensions to the generated
$C^*$-algebra \cite[Definition~2.1.1]{Arveson}.

Stinespring's theorem is equivalent to Naimark's dilation theorem
for commutative domains \cite{Naimark} (see \cite[pp.~211--213]{Stinespring}),
and, in general, to Sz.-Nagy's general dilation theorem
\cite{Szafraniec1983} (see \cite[pp.~289--290]{Szafraniec2010}).

A linear preserver problem asks which linear maps preserve a specified
function, set, or relation on an algebra. The aim is to recover the
structure of the map from the data it preserves.

Preservation of spectral data is a classical example. Marcus and
Moyls~\cite[Theorem~3]{MM} characterize linear maps on $M_d(\C)$
preserving eigenvalues, counted with multiplicities, as similarities or
similarities composed with transposition. Their paper also treats rank
and determinant preservation.
Preservation of zero products is another extensively studied condition;
see, for example, \cite{CKLW,ABEV}.

These problems have matrix-level counterparts, in which the preservation
requirement is imposed on every entrywise amplification of the map.
For a complex-linear map $\Phi\colon \mathscr{A}\to \mathscr{B}$, write
\[
\Phi^{(n)}\colon M_n(\mathscr{A})\longrightarrow M_n(\mathscr{B}),\qquad
\Phi^{(n)}([a_{ij}])=[\Phi(a_{ij})].
\]
We say that $\Phi$ \emph{completely preserves} an operator property if every
$\Phi^{(n)}$ preserves it. The condition concerns the whole matrix, not its
entries.

Examples from the literature include complete preservation of spectra,
invertibility, and orthogonality. For surjective maps between standard
operator algebras on complex Banach spaces, Hou and Huang~\cite{HH}
characterize algebra isomorphisms by complete spectrum preservation.
They also characterize ring isomorphisms by unitality and preservation
of invertibility in both directions at every matrix level.
For bounded linear maps between $C^*$-algebras, Garc\'es~\cite[Theorem~14]{Garces}
characterizes complete orthogonality preservers as weighted TRO
homomorphisms and shows that the second matrix level already suffices.
These results illustrate how matrix-level preservation conditions
determine the algebraic structure of a map.

\section{Main result}\label{sec:main-results}

From the viewpoint of complete preserver problems, complete positivity
is exactly complete preservation of positivity:
\[
x\in M_n(\mathscr{A}),\quad x\geq0
\quad\Longrightarrow\quad \Phi^{(n)}(x)\geq0
\qquad(n\geq1).
\]
For maps with values in $\B(\Hh)$, Stinespring's theorem
(see Theorem~\ref{thm:stinespring}) gives a representation-theoretic
solution to the complete preserver problem for positivity.

We ask which maps admit a corresponding representation when subnormality
takes the place of positivity.

An operator is subnormal if it is the restriction of a normal operator to
an invariant subspace; see Section~\ref{sec:prerequisites} for the intrinsic
$C^*$-algebra formulation. For a bounded linear map
$\Phi\colon \mathscr{A}\to\B(\Hh)$, let
\begin{equation}\label{eq:canonical-extension}
\Phihat\colon \mathscr{A}^{**}\longrightarrow\B(\Hh)
\end{equation}
denote its canonical weak-$*$ continuous extension and put $H=\Phihat(1)$.
If $\mathscr{A}$ is unital, then $H=\Phi(1)$.

The following result, like Stinespring's theorem, characterizes the maps
by a representation on a larger Hilbert space.

\begin{theorem}
\label{thm:subnormal-stinespring}
Let $\mathscr{A}$ be a $C^*$-algebra, let $\Hh$ be a Hilbert space,
and let $\Phi\colon \mathscr{A}\to\B(\Hh)$ be a bounded linear map.
Set $H=\Phihat(1)$, where $\Phihat$ is the extension
in \eqref{eq:canonical-extension}.
The following conditions are equivalent:
\begin{enumerate}
\item
For every $n\geq1$, the amplification $\Phi^{(n)}$ maps
normal elements of $M_n(\mathscr{A})$ to subnormal operators.

\item
For every $n\geq1$, the amplification $\Phi^{(n)}$
preserves subnormality.

\item
The operator $H$ is subnormal, and its minimal normal
extension $N\in\B(\Kk)$, where $\Kk\supseteq\Hh$, admits
a $*$-homomorphism $\rho\colon \mathscr{A}\to\B(\Kk)$ such that, for every $a\in \mathscr{A}$,
\begin{equation}\label{eq:extension-stinespring}
[N,\rho(a)]=0,\qquad N\rho(a)\Hh\subseteq\Hh,\qquad
\Phi(a)=\bigl(N\rho(a)\bigr)|_{\Hh}.
\end{equation}
\end{enumerate}

Moreover, let $P$ be the orthogonal projection onto
$\overline{\ran N}$. For the fixed minimal normal extension $N$,
there is a unique choice of $\rho$ satisfying
\[
\rho(\mathscr{A})\subseteq P\B(\Kk)P.
\]
This representation is nondegenerate on $P\Kk$, and
\begin{equation}\label{eq:subnorm-stinespring}
\|\Phi\|_{\mathrm{cb}}=\|\Phi\|=\|H\|.
\end{equation}
If $\mathscr{A}$ is unital, the boundedness hypothesis may be omitted,
$H=\Phi(1)$, and the canonical representation satisfies $\rho(1)=P$.
\end{theorem}

Theorem~\ref{thm:subnormal-stinespring} gives a related representation
principle for complete preservation of subnormality. The map $\Phi$ need
not be positive; it is realized as an invariant restriction of a
representation multiplied by a commuting normal operator.
Specifically, \eqref{eq:extension-stinespring} gives
$N\rho(a)\xi=\Phi(a)\xi$ for $\xi\in\Hh$, so that $\Hh$ is a common
invariant subspace for the operators $N\rho(a)$.
Neither $\rho(a)$ nor $N^*$ is required to leave $\Hh$ invariant.
The same space $\Kk$ also supports joint normal extensions for images of
finite commuting families of normal elements.

\subsection*{AI Disclosure}
The authors used ChatGPT while preparing this manuscript.

\section{Prerequisites}\label{sec:prerequisites}

\subsection{Notation and operator terminology}
All Hilbert spaces are complex, and inner products are linear in the
second variable. All linear maps between $C^*$-algebras are understood
to be complex-linear. We use uppercase Latin letters for operators
on Hilbert spaces and lowercase Latin letters for elements of
abstract $C^*$-algebras. We write $\C$ for the complex numbers and
$\mathbb N_0=\{0,1,2,\ldots\}$.

The algebra of bounded operators on a Hilbert space $\Hh$ is denoted
by $\B(\Hh)$, with identity $I$. For $T\in\B(\Hh)$, the symbols
$T^*$, $\ker T$, and $\ran T$ denote its adjoint, kernel, and range.
The orthogonal projection onto a closed
subspace $\Ee$ is denoted by $P_{\Ee}$. We identify $M_n(\B(\Hh))$
with $\B(\Hh^n)$, where $\Hh^n$ is the direct sum of $n$ copies of
$\Hh$, and write $I_n$ for the identity matrix in $M_n(\C)$.
An operator $T$ is hyponormal if $T^*T-TT^*\geq0$, and is a contraction
if $\|T\|\leq1$.

We write $[S,T]=ST-TS$. Operators $S,T$ doubly commute if
$ST=TS$ and $ST^*=T^*S$. For a set $\mathcal F\subseteq\B(\Hh)$,
its commutant is
$\mathcal F'=\{X\in\B(\Hh):XF=FX\text{ for all }F\in\mathcal F\}$.
The notation $C^*(\mathcal F)$ denotes the $C^*$-algebra generated by
$\mathcal F$, without adjoining an ambient identity.
A $*$-homomorphism $\pi\colon \mathscr{A}\to\B(\Kk)$ is nondegenerate if
$\overline{\operatorname{span}}\pi(\mathscr{A})\Kk=\Kk$.

An operator commuting with a positive element commutes with its
continuous functions and spectral projections. We also use norm
approximation of self-adjoint elements of a von Neumann algebra by
finite linear combinations of mutually orthogonal spectral projections;
see \cite{Blackadar,Conway}.

For a bounded linear map $\Phi\colon \mathscr{A}\to \mathscr{B}$, its completely bounded
norm is $\|\Phi\|_{\mathrm{cb}}=\sup_{n\geq1}\|\Phi^{(n)}\|$;
the map is completely bounded when this supremum is finite.
Elements $a,b$ of a $C^*$-algebra are range-orthogonal if $a^*b=0$.

\subsection{Subnormal operators and moment criteria}
For a subnormal operator $T\in\B(\Hh)$, its minimal normal extension
$N\in\B(\Kk)$, $\Kk\supseteq\Hh$, is characterized by
\begin{equation}\label{eq:minimal}
\Kk=\overline{\operatorname{span}}
\{N^{*j}\xi:j\geq0,\ \xi\in\Hh\}.
\end{equation}
It is unique up to a unitary fixing $\Hh$ and has norm $\|T\|$;
see \cite{Bram,Conway}.

We will use the following form of Bram's commutant theorem; see
\cite[Theorem~8, p.~89]{Bram} and
\cite[Chapter~II, Theorem~10.5]{Conway}.
\begin{theorem}\label{thm:bram-commutant}
Let $T\in\B(\Hh)$ be subnormal, let $N\in\B(\Kk)$ be its minimal
normal extension, and put $P=P_{\Hh}$. Then restriction to $\Hh$ is a
unital $*$-isomorphism
\[
C^*(N)'\cap\{P\}'\longrightarrow C^*(T)',\qquad X\longmapsto X|_{\Hh}.
\]
\end{theorem}

For $m\geq1$, define $\mathcal B_m\colon\B(\Hh)\to M_{m+1}(\B(\Hh))$ by
$\mathcal B_m(T)=[T^{*j}T^i]_{i,j=0}^m$. The Bram--Halmos criterion gives
\begin{equation}\label{eq:BH}
T\text{ is subnormal}
\quad\Longleftrightarrow\quad
\mathcal B_m(T)\geq0
\quad(m\geq1).
\end{equation}
In particular,
$\mathcal B_1(T)=\left(\begin{smallmatrix}I&T^*\\T&T^*T\end{smallmatrix}\right)$,
so subnormality implies hyponormality.
For an element of a nonunital $C^*$-algebra the criterion is interpreted
in its unitization. It follows that subnormality is independent of the
faithful representation, is preserved by $*$-homomorphisms, and is closed
under norm-bounded strong-$*$ limits. For a nonunital homomorphism the
zero summand causes no difficulty: its moment matrix is positive as well.
These consequences follow by applying the homomorphism, or taking limits,
in each fixed matrix in \eqref{eq:BH}. They also show that normality and
subnormality agree for matrices over $\C$, since a positive
self-commutator in finite dimension has trace zero.
Preservation of subnormality at any matrix level may therefore be tested
after a faithful representation of the target algebra.

\begin{lemma}\label{lem:products}
If $R,T\in\B(\Hh)$ are subnormal and doubly commute, then $RT$ is
subnormal. In particular, this holds when $R$ is normal and commutes with
$T$ and $T^*$.
\end{lemma}
\begin{proof}
By \cite[Theorem~3]{Ando}, $R$ and $T$ admit commuting normal
extensions on a common Hilbert space. Their product is a normal
extension of $RT$.
\end{proof}

\subsection{Biduals and operator topologies}
For a $C^*$-algebra $\mathscr{A}$, we identify $\mathscr{A}$ with its canonical image in
the enveloping von Neumann algebra $\mathscr{A}^{**}$. A $*$-homomorphism need
not be unital. For a von Neumann algebra $\mathscr{M}$, its predual and the
positive cone of its predual are denoted by $\mathscr{M}_*$ and $\mathscr{M}_*^+$.

The second adjoint of a bounded linear map $T\colon \mathscr{A}\to \mathscr{B}$ is
$T^{**}\colon \mathscr{A}^{**}\to \mathscr{B}^{**}$. If $\mathscr{B}=\B(\Hh)$,
the canonical extension $\widehat T\colon\mathscr{A}^{**}\to\B(\Hh)$
of $T$, as in \eqref{eq:canonical-extension}, is characterized by
\[
\omega(\widehat T(x))=\langle x,\omega\circ T\rangle
\qquad(x\in \mathscr{A}^{**},\ \omega\in\B(\Hh)_*).
\]
Moreover, $\|\widehat T\|=\|T\|$.

A net $(T_\alpha)$ in $\B(\Hh)$ converges strongly to $T$ if
$T_\alpha\xi\to T\xi$ for every $\xi\in\Hh$; it converges strongly-$*$
if the same also holds for the adjoints. We write
$\operatorname*{s\!-\!lim}$ for limits in the strong operator topology.

On a von Neumann algebra $\mathscr{M}$, the $\sst$ topology is generated by
\[
x\longmapsto\bigl(\omega(x^*x+xx^*)\bigr)^{1/2},\qquad\omega\in \mathscr{M}_*^+.
\]
On norm-bounded sets it agrees with strong-$*$ convergence in a faithful
normal representation. Multiplication is jointly continuous there.

The following continuity fact is a standard consequence of the
Mackey-topology description of bounded sets in a von Neumann algebra.
We give the short duality argument to make the exact topology explicit.

\begin{lemma}\label{lem:continuity}
For every bounded linear map $T\colon \mathscr{A}\to \mathscr{B}$, its second adjoint
$T^{**}\colon \mathscr{A}^{**}\to \mathscr{B}^{**}$ is $\sst$ continuous on norm-bounded
sets. If $\mathscr{B}=\B(\Hh)$, the same holds for the canonical normal extension
$\widehat T\colon \mathscr{A}^{**}\to\B(\Hh)$.
\end{lemma}
\begin{proof}
Let $\Theta\colon\mathscr{M}\to\mathscr{N}$ be a bounded weak-$*$
continuous map between von Neumann algebras. Its preadjoint
$\Theta_*\colon\mathscr{N}_*\to\mathscr{M}_*$ carries weakly compact
sets to weakly compact sets. Thus $\Theta$ is continuous for the
corresponding Mackey topologies. By \cite[p.~93, Corollary~1]{Aarnes},
these agree with $\sst$ on norm-bounded sets. Apply this observation to $T^{**}$
and, when $\mathscr{B}=\B(\Hh)$, to the canonical weak-$*$ continuous
extension $\widehat T$.
\end{proof}

To pass preservation properties to the bidual, we use the following
version of Kaplansky's density theorem for normal elements. The unital
case is contained in Hadwin's simultaneous approximation lemma
\cite[Lemma~4]{Hadwin}; we include the reduction to nonunital algebras.

\begin{lemma}\label{lem:density}
Let $\mathscr{A}$ be a $C^*$-algebra. Every normal contraction in $\mathscr{A}^{**}$ is
the $\sst$ limit of a net of normal contractions in $\mathscr{A}$.
\end{lemma}
\begin{proof}
If $\mathscr{A}$ is unital, the assertion follows from \cite[Lemma~4]{Hadwin}
in the universal representation, since the strong-$*$ and $\sst$
topologies agree on bounded sets.

Suppose that $\mathscr{A}$ is nonunital and let $x\in \mathscr{A}^{**}$ be a normal
contraction. Let $\widetilde{\mathscr{A}}$ be the unitization of $\mathscr{A}$ and identify
$\widetilde{\mathscr{A}}^{**}$ with $\mathscr{A}^{**}\oplus\C$, where
$a+\lambda1_{\widetilde{\mathscr{A}}}$ corresponds to
$(a+\lambda1_{\mathscr{A}^{**}},\lambda)$. By the unital case, there is a net
of normal contractions
$b_\alpha=a_\alpha+\lambda_\alpha1_{\widetilde{\mathscr{A}}}\in\widetilde{\mathscr{A}}$
converging $\sst$ to $(x,0)$. Thus $\lambda_\alpha\to0$ and
$a_\alpha\to x$ in $\sst$. Since each $a_\alpha$ is normal and
$\|a_\alpha\|\leq1+|\lambda_\alpha|$, the elements
$c_\alpha=(1+|\lambda_\alpha|)^{-1}a_\alpha$ are normal contractions
in $\mathscr{A}$ converging $\sst$ to $x$.
\end{proof}

\section{Automatic boundedness and passage to the bidual}
\label{sec:boundedness-bidual}

\begin{lemma}\label{lem:auto}
Let $\mathscr{A},\mathscr{B}$ be $C^*$-algebras with $\mathscr{A}$ unital, and let
$\Phi\colon \mathscr{A}\to \mathscr{B}$ be complex-linear. Assume that $\Phi^{(2)}(u)$
is hyponormal for every unitary $u\in M_2(\mathscr{A})$, and put $h=\Phi(1)$.
Then $\Phi$ is bounded and $\|\Phi\|=\|h\|$.
If $\Phi^{(n)}$ maps unitaries to hyponormal elements for every $n\geq1$,
then $\Phi$ is completely bounded and $\|\Phi\|_{\mathrm{cb}}=\|h\|$.
\end{lemma}
\begin{proof}
For a unitary $u\in \mathscr{A}$, apply the hypothesis to
\[
w_u=\begin{pmatrix}0&u\\1&0\end{pmatrix}.
\]
Hyponormality of $\Phi^{(2)}(w_u)$ gives
\[
\Phi(u)\Phi(u)^*\leq h^*h,\qquad
\Phi(u)^*\Phi(u)\geq hh^*.
\]
Thus $\|\Phi(u)\|\leq\|h\|$. Every $a\in \mathscr{A}$ is a linear combination
of four unitaries with sum of absolute coefficients at most $2\|a\|$,
so $\Phi$ is bounded with $\|\Phi\|\leq2\|h\|$.
The Russo--Dye theorem \cite{RD} now yields
\[
\|\Phi\|=\sup_{u\text{ unitary}}\|\Phi(u)\|=\|h\|.
\]

If every amplification has the stated preservation property,
apply the first part to $\Phi^{(n)}$, using the hypothesis at level $2n$.
Then
\[
\|\Phi^{(n)}\|=\|I_n\otimes h\|=\|h\|\qquad(n\geq1),
\]
and hence $\|\Phi\|_{\mathrm{cb}}=\|h\|$.
\end{proof}

Lemmas~\ref{lem:continuity} and \ref{lem:density} yield the following
extension of preservation properties to the bidual.

\begin{proposition}\label{prop:bidual}
Let $\mathscr{A}$ and $\mathscr{B}$ be $C^*$-algebras, let $T\colon \mathscr{A}\to \mathscr{B}$ be a bounded
linear map, and fix $n\geq1$. Then $T^{(n)}$ maps normal elements to
subnormal elements if and only if $(T^{**})^{(n)}$ does.
If $\mathscr{B}=\B(\Hh)$, the same equivalence holds for the canonical weak-$*$
continuous extension $\widehat T\colon \mathscr{A}^{**}\to\B(\Hh)$.
\end{proposition}
\begin{proof}
Since $T^{(n)}$ is bounded and, under the canonical identifications,
$M_n(\mathscr{A})^{**}=M_n(\mathscr{A}^{**})$ and
$(T^{(n)})^{**}=(T^{**})^{(n)}$, it suffices to consider $n=1$.
Suppose that $T$ maps normal elements to subnormal elements,
and let $x\in \mathscr{A}^{**}$ be a normal contraction. By
Lemma~\ref{lem:density}, there is a net $(a_\alpha)$ of normal
contractions in $\mathscr{A}$ converging $\sst$ to $x$.
Lemma~\ref{lem:continuity} gives $T(a_\alpha)\to T^{**}(x)$ in $\sst$.
Since $(T(a_\alpha))$ is norm bounded, the positivity of each moment
matrix in \eqref{eq:BH} passes to the strong-$*$ limit. Thus $T^{**}(x)$
is subnormal.
Scaling gives the assertion for every normal $x\in \mathscr{A}^{**}$.

The reverse implication follows by restriction to $\mathscr{A}$.
For $\mathscr{B}=\B(\Hh)$, the same argument applies to $\widehat T$ and its
amplifications, using the concrete form of
Lemma~\ref{lem:continuity}.
\end{proof}

\section{Proof of Theorem~\ref{thm:subnormal-stinespring}}\label{sec:proof-stinespring}

We first isolate a finite-dimensional test for two arbitrary subnormal
operators. The test exploits the initial moments of two cycles before
mixed paths return to their common vertex.

For $k\geq1$, define $Q_k\colon\B(\Hh)\to\B(\Hh)$ by
\begin{equation}\label{eq:Qdef}
Q_k(T)=\sum_{j=0}^k(-1)^j\binom{k}{j}T^{*j}T^j.
\end{equation}
By Agler's criterion \cite[Theorem~3.1]{Agler}, a contraction $T$
is subnormal if and only if $Q_k(T)\geq0$ for every $k\geq1$.
This characterization is based on Embry's moment criterion
\cite{Embry}. Since $Q_{k+1}(T)=Q_k(T)-T^*Q_k(T)T$, every subnormal
contraction satisfies
\begin{equation}\label{eq:Qbounds}
0\leq Q_k(T)\leq Q_1(T)=I-T^*T\qquad(k\geq1).
\end{equation}

\begin{lemma}\label{lem:binomial}
Let $T\in\B(\Hh)$ be a subnormal contraction. Then, for every $k\geq1$,
\begin{equation}\label{eq:TQT}
T^*Q_k(T)T\leq\frac{k^k}{(k+1)^{k+1}}I.
\end{equation}
\end{lemma}
\begin{proof}
Let $N\in\B(\Kk)$ be the minimal normal extension of $T$, so that
$\|N\|=\|T\|\leq1$. Invariance of $\Hh$ gives
\[
T^*Q_k(T)T
=P_{\Hh}N^*N(I-N^*N)^k|_{\Hh}.
\]
The assertion follows by functional calculus from
$\max_{0\leq t\leq1}t(1-t)^k=k^k/(k+1)^{k+1}$.
Equality holds for $T=(k+1)^{-1/2}I$.
\end{proof}

The next theorem derives zero-product preservation from matricial
subnormality. Its proof uses two cyclic shifts sharing one coordinate
and a quantitative estimate obtained from the binomial moment
conditions.

\begin{theorem}\label{thm:zero-products}
Let $S,T\in\B(\Hh)$. If $X\otimes S+Y\otimes T$ is subnormal
for every $n\geq1$ and every pair of normal matrices
$X,Y\in M_n(\C)$, then $ST=TS=0$.
\end{theorem}
\begin{proof}
The scalar choices $(1,0)$ and $(0,1)$ show that $S$ and $T$ are
subnormal. If $S=T=0$, the conclusion is immediate. Fix $m\geq2$
and choose an orthonormal basis
$e_0,e_1,\ldots,e_{m-1},f_1,\ldots,f_{m-1}$ of $\C^{2m-1}$.
Let $C_m$ and $D_m$ be the cyclic shifts, in the indicated orders, on
$\operatorname{span}\{e_0,\ldots,e_{m-1}\}$ and
$\operatorname{span}\{e_0,f_1,\ldots,f_{m-1}\}$, respectively,
extended by zero on the corresponding orthogonal complements. Set
$Z_m=C_m\otimes S+D_m\otimes T$. The hypothesis implies that $Z_m$
is subnormal. We prove that
\begin{equation}\label{eq:cycle-bound}
\max\{\|ST\|^2,\|TS\|^2\}
\leq\|S^*S+T^*T\|^2\frac{(m-1)^{m-1}}{m^m}.
\end{equation}

The matrices $C_m$ and $D_m$ are normal contractions. By homogeneity,
we may assume that $\|S^*S+T^*T\|=1$. Relative to suitable orthogonal
decompositions of the domain and range, $Z_m$ is the direct sum of
the blocks $\binom{S}{T}$, $(S\ T)$, and $m-2$ copies of each of
$S$ and $T$. Since $S$ and $T$ are hyponormal,
\[
\|Z_m\|^2
=\max\{\|S^*S+T^*T\|,\|SS^*+TT^*\|\}=1.
\]
Thus $S$, $T$, and $Z_m$ are subnormal contractions.

Define $J\colon\Hh\to\C^{2m-1}\otimes\Hh$ by
$J\xi=e_0\otimes\xi$. Neither cycle returns to the common
coordinate before step $m$, and hence
\[
Z_m^jJ\xi=e_j\otimes S^j\xi+f_j\otimes T^j\xi
\qquad(1\leq j<m,\ \xi\in\Hh).
\]
For $1\leq k<m$, this identity and the positivity of $Q_k(Z_m)$ give
\[
0\leq J^*Q_k(Z_m)J=Q_k(S)+Q_k(T)-I.
\]
Together with~\eqref{eq:Qbounds}, this yields
$S^*S\leq I-Q_k(S)\leq Q_k(T)$. Lemma~\ref{lem:binomial}
therefore gives
\[
(ST)^*(ST)=T^*S^*ST
\leq T^*Q_k(T)T
\leq\frac{k^k}{(k+1)^{k+1}}I.
\]
Taking $k=m-1$ and interchanging $S$ and $T$ proves
\eqref{eq:cycle-bound} under the normalization. Undoing the
normalization proves it in general. Since
$(m-1)^{m-1}/m^m\leq1/m$ and
$\|S^*S+T^*T\|\leq\|S\|^2+\|T\|^2$, the right-hand side of
\eqref{eq:cycle-bound} is at most $(\|S\|+\|T\|)^4/m$.
Letting $m\to\infty$ in~\eqref{eq:cycle-bound} gives $ST=TS=0$.
\end{proof}

\begin{proposition}
\label{prop:algebraic-subnormality}
Let $\Phi\colon \mathscr{A}\to\B(\Hh)$ be a linear map, where
$\mathscr{A}$ is a unital $C^*$-algebra and $\Hh$ is a Hilbert space.
Set $H=\Phi(1)$. The following are equivalent.
\begin{enumerate}
\item
For every $n\geq1$ and every normal $x\in M_n(\mathscr{A})$,
the operator $\Phi^{(n)}(x)$ is subnormal.

\item
$\Phi$ completely preserves subnormality.

\item
The operator $H$ is subnormal and the identities
\begin{align}
\Phi(a)\Phi(b)&=H\Phi(ab)=\Phi(ab)H,\label{eq:alg-mult-stinespring}\\
\Phi(a)^*\Phi(b)&=H^*\Phi(a^*b)\label{eq:alg-star-stinespring}
\end{align}
hold for all $a,b\in\mathscr{A}$.
\end{enumerate}
If these conditions hold, then $\Phi$ is bounded and preserves zero products
and range-orthogonality.
\end{proposition}

For bounded maps, zero-product preservation gives \eqref{eq:alg-mult-stinespring}
by \cite[Theorem~2.3(1)]{GT}, and range-orthogonality preservation gives
\eqref{eq:alg-star-stinespring} by Schweizer's theorem
\cite[Theorem~3.6]{LCLW}. Here the two-cycle test yields both identities
from condition~\textup{(i)}.

We first establish the extension step used for
\textup{(iii)}$\Rightarrow$\textup{(ii)}.

\begin{lemma}\label{lem:algebraicmodel}
Let $\Phi\colon\mathscr{A}\to\B(\Hh)$ be a bounded linear map, where
$\mathscr{A}$ is a unital $C^*$-algebra and $\Hh$ is a Hilbert space.
Suppose that $H=\Phi(1)$ is subnormal. Let $N\in\B(\Kk)$ be its minimal
normal extension and $P$ the orthogonal projection onto $\overline{\ran N}$.
Then \eqref{eq:alg-mult-stinespring}--\eqref{eq:alg-star-stinespring}
hold for all $a,b\in\mathscr{A}$ if and only if there is a
$*$-homomorphism $\rho\colon\mathscr{A}\to P\B(\Kk)P$ satisfying
\eqref{eq:extension-stinespring} for every $a\in\mathscr{A}$.
In this case $\rho$ is unique, $\rho(1)=P$, and $\Phi$ completely
preserves subnormality.
\end{lemma}
\begin{proof}
If $H=0$, then $N=0$ by minimality. Both conditions reduce to $\Phi=0$,
and the conclusion holds with $\rho=0$. Assume henceforth that $H\ne0$.

Suppose that \eqref{eq:alg-mult-stinespring}--\eqref{eq:alg-star-stinespring}
hold. Define $\Lambda\colon\mathscr{A}\to\B(\Hh)$ by
$\Lambda(a)=H^*\Phi(a)$ and set
\[
\Ee=\overline{\operatorname{span}}\{\Phi(b)\xi:b\in\mathscr{A},\ \xi\in\Hh\}.
\]
By \eqref{eq:alg-star-stinespring}, the matrices
$[\Lambda(a_i^*a_j)]_{i,j}=[\Phi(a_i)^*\Phi(a_j)]_{i,j}$ are positive.
Hence $\Lambda$ is completely positive. In the minimal Stinespring
construction \cite{Stinespring} (see also \cite[Theorem~1.1.1]{Arveson}),
the linear map from the algebraic tensor product $\mathscr{A}\odot\Hh$
to $\Ee$ given by $b\otimes\xi\mapsto\Phi(b)\xi$ induces a unitary from
the Stinespring space onto $\Ee$. The resulting unital $*$-representation
$\pi\colon\mathscr{A}\to\B(\Ee)$ satisfies
\[
\pi(a)\Phi(b)=\Phi(ab)\qquad(a,b\in\mathscr{A}).
\]

By \eqref{eq:alg-mult-stinespring}, $H$ commutes with $\Phi(\mathscr{A})$,
so $\Ee$ is invariant under $H$. For $T_0=H|_{\Ee}\in\B(\Ee)$ we have
\[
T_0\pi(a)\Phi(b)=H\Phi(ab)=\Phi(ab)H=\pi(a)T_0\Phi(b).
\]
The ranges of the operators $\Phi(b)$ span a dense subspace of $\Ee$.
Thus $T_0$ commutes with $\pi(\mathscr{A})$, and taking adjoints gives
$\pi(\mathscr{A})\subseteq C^*(T_0)'$.

Set
\[
\mathcal L=\overline{\operatorname{span}}
\{N^{*j}\eta:j\geq0,\ \eta\in\Ee\}.
\]
The space $\mathcal L$ is invariant under $N^*$ by construction and under
$N$ because $NN^{*j}\eta=N^{*j}T_0\eta$ for $\eta\in\Ee$. Thus
$N|_{\mathcal L}$ is the minimal normal extension of $T_0$.
Bram's commutant theorem (see Theorem~\ref{thm:bram-commutant}) lifts
$\pi$ to a unital $*$-representation $\mathscr{A}\to\B(\mathcal L)$
whose range commutes with $N|_{\mathcal L}$. Extending this lift by zero
on $\mathcal L^\perp$ gives a $*$-representation
$\rho\colon\mathscr{A}\to\B(\Kk)$ whose range commutes with $N$. Hence
\[
N\rho(a)\xi=\rho(a)H\xi=\pi(a)H\xi=\Phi(a)\xi
\qquad(a\in\mathscr{A},\ \xi\in\Hh).
\]
In particular, $\Ee\subseteq P\Kk$. Since $P\Kk$ reduces $N$, we have
$\mathcal L\subseteq P\Kk$. Applying $N$ to the dense span in
\eqref{eq:minimal} and taking closures gives
\begin{equation}\label{eq:supportdense}
P\Kk=\overline{\operatorname{span}}
\{N^{*j}H\xi:j\geq0,\ \xi\in\Hh\}.
\end{equation}
Since $H\Hh\subseteq\Ee$, the right-hand side is contained in
$\mathcal L$. Thus $\mathcal L=P\Kk$. The extension by zero therefore
allows us to regard $\rho$ as a $*$-homomorphism into $P\B(\Kk)P$,
with $\rho(1)=P$.

Let $\widetilde\rho\colon\mathscr{A}\to P\B(\Kk)P$ be any
$*$-homomorphism satisfying \eqref{eq:extension-stinespring}. The ranges
of $\rho$ and $\widetilde\rho$ commute with $N$ and, by self-adjointness,
with $N^*$. For $a\in\mathscr{A}$, $j\geq0$, and $\xi\in\Hh$, we have
\begin{equation}\label{eq:uniqueness}
\rho(a)N^{*j}H\xi=\widetilde\rho(a)N^{*j}H\xi=N^{*j}\Phi(a)\xi.
\end{equation}
By \eqref{eq:supportdense}, $\rho(a)$ and $\widetilde\rho(a)$ agree on
$P\Kk$. Both vanish on $(P\Kk)^\perp$, so $\widetilde\rho=\rho$.

Conversely, suppose that a $*$-homomorphism
$\rho\colon\mathscr{A}\to P\B(\Kk)P$ satisfies \eqref{eq:extension-stinespring}.
Its range commutes with $N$ and $N^*$. Since $\Hh$ is invariant under $N$
and every $N\rho(a)$, we have
\begin{align*}
\Phi(a)\Phi(b)&=(N^2\rho(ab))|_{\Hh}=H\Phi(ab)=\Phi(ab)H,\\
\Phi(a)^*\Phi(b)&=P_{\Hh}N^*N\rho(a^*b)|_{\Hh}=H^*\Phi(a^*b)
\end{align*}
for all $a,b\in\mathscr{A}$.

If $x\in M_n(\mathscr{A})$ is subnormal, then $\rho^{(n)}(x)$ is
subnormal by \eqref{eq:BH} and doubly commutes with the normal operator
$I_n\otimes N$. Lemma~\ref{lem:products} gives subnormality of
$(I_n\otimes N)\rho^{(n)}(x)$. By \eqref{eq:extension-stinespring},
$\Phi^{(n)}(x)$ is an invariant restriction of this product, hence is
subnormal.
\end{proof}

\begin{proof}[Proof of Proposition~\ref{prop:algebraic-subnormality}]
The implication \textup{(ii)}$\Rightarrow$\textup{(i)} is immediate.
Assume that \textup{(i)} holds. Since unitaries are normal and subnormal
operators are hyponormal, Lemma~\ref{lem:auto} gives the boundedness of
$\Phi$. Let $\Psi=\Phihat\colon\mathscr{A}^{**}\to\B(\Hh)$ be its
canonical extension.
By Proposition~\ref{prop:bidual}, every matrix amplification of $\Psi$
maps normal elements to subnormal operators. In particular, $H=\Psi(1)$
is subnormal.

Let $e,f\in M_2(\mathscr{A}^{**})$ be orthogonal projections and put
$S=\Psi^{(2)}(e)$ and $T=\Psi^{(2)}(f)$. For every $n\geq1$
and all normal $X,Y\in M_n(\C)$, the element $X\otimes e+Y\otimes f$ is
normal, so its image $X\otimes S+Y\otimes T$ under $\Psi^{(2n)}$ is
subnormal.
Theorem~\ref{thm:zero-products} gives $ST=TS=0$. Since $S$ and $T$ are
hyponormal, $\ker S\subseteq\ker S^*$ and $\ker T\subseteq\ker T^*$.
Thus $S^*T=T^*S=0$.

Put $H_2=I_2\otimes H$. Taking $f=1-e$ and using
$\Psi^{(2)}(1-e)=H_2-S$ gives
\[
S^2=H_2S=SH_2,\qquad S^*S=H_2^*S.
\]
For a self-adjoint $x\in M_2(\mathscr{A}^{**})$ of finite spectrum, write
$x=\sum_j\lambda_je_j$ in terms of its spectral projections. The mixed
products vanish, so
\begin{align*}
\Psi^{(2)}(x)^2&=H_2\Psi^{(2)}(x^2)=\Psi^{(2)}(x^2)H_2,\\
\Psi^{(2)}(x)^*\Psi^{(2)}(x)&=H_2^*\Psi^{(2)}(x^2).
\end{align*}
By norm approximation, these identities hold for every self-adjoint
$x\in M_2(\mathscr{A}^{**})$.

For $a,b\in\mathscr{A}^{**}$ with $a=a^*$, take
\[
x=\begin{pmatrix}a&b\\ b^*&0\end{pmatrix}.
\]
The $(1,2)$ entries of the two identities give
\begin{align*}
\Psi(a)\Psi(b)&=H\Psi(ab)=\Psi(ab)H,\\
\Psi(a)^*\Psi(b)&=H^*\Psi(a^*b).
\end{align*}
The first identity is linear in $a$ and the second is conjugate-linear
in $a$, so both hold for arbitrary $a,b\in\mathscr{A}^{**}$. Restricting to
$\mathscr{A}$ gives \textup{(iii)}.

Assume now that \textup{(iii)} holds and define
$\Lambda\colon\mathscr{A}\to\B(\Hh)$ by $\Lambda(a)=H^*\Phi(a)$.
By \eqref{eq:alg-star-stinespring}, $\Lambda(a^*a)=\Phi(a)^*\Phi(a)\geq0$,
so $\Lambda$ is positive. Therefore
\[
\Phi(a)^*\Phi(a)=\Lambda(a^*a)
\leq\|a\|^2\Lambda(1)=\|a\|^2H^*H
\qquad(a\in\mathscr{A}).
\]
Thus $\Phi$ is bounded, and Lemma~\ref{lem:algebraicmodel} gives
\textup{(ii)}. The identities \eqref{eq:alg-mult-stinespring} and
\eqref{eq:alg-star-stinespring} give preservation of zero products and
range-orthogonality, respectively.
\end{proof}

\begin{proof}[Proof of Theorem~\ref{thm:subnormal-stinespring}]
The implication \textup{(ii)}$\Rightarrow$\textup{(i)} is immediate.
Assume \textup{(i)}. If $\mathscr{A}$ is unital, the images of unitaries
under $\Phi^{(2)}$ are hyponormal, so Lemma~\ref{lem:auto} gives
boundedness. Proposition~\ref{prop:bidual} extends \textup{(i)} to
$\Phihat$. Proposition~\ref{prop:algebraic-subnormality} applied to
$\Phihat$ gives subnormality of $H$ and
\eqref{eq:alg-mult-stinespring}--\eqref{eq:alg-star-stinespring} for
$\Phihat$.

Let $N\in\B(\Kk)$ be the minimal normal extension of $H$ and $P$ the
orthogonal projection onto $\overline{\ran N}$.
Lemma~\ref{lem:algebraicmodel} gives a $*$-homomorphism
$\widetilde\rho\colon\mathscr{A}^{**}\to P\B(\Kk)P$ satisfying
\eqref{eq:extension-stinespring} with $\Phihat$ in place of $\Phi$.
Its restriction $\rho\colon\mathscr{A}\to P\B(\Kk)P$ proves
\textup{(iii)}.

Let $\sigma\colon\mathscr{A}\to P\B(\Kk)P$ be any $*$-homomorphism
satisfying \eqref{eq:extension-stinespring}. The ranges of $\sigma$
and $\rho$ commute with $N$ and, by self-adjointness, with $N^*$.
As in \eqref{eq:uniqueness}, for $a\in\mathscr{A}$, $j\geq0$, and
$\xi\in\Hh$ we obtain
\[
\sigma(a)N^{*j}H\xi=N^{*j}\Phi(a)\xi=\rho(a)N^{*j}H\xi.
\]
By \eqref{eq:supportdense}, the two operators agree on $P\Kk$.
They vanish on $(P\Kk)^\perp$, so $\sigma=\rho$.

Let $(e_\lambda)$ be a positive contractive approximate identity of
$\mathscr{A}$. Since $e_\lambda\to1$ in the $\sst$ topology of
$\mathscr{A}^{**}$, Lemma~\ref{lem:continuity} gives
$\Phi(e_\lambda)\to H$ strongly. The preceding identity gives
$\rho(e_\lambda)N^{*j}H\xi\to N^{*j}H\xi$ for $j\geq0$ and $\xi\in\Hh$.
By \eqref{eq:supportdense} and contractivity, $\rho(e_\lambda)\to P$
strongly. Thus $\rho$ is nondegenerate on $P\Kk$. If $\mathscr{A}$ is
unital, the constant approximate identity gives $\rho(1)=P$.

Assume now \textup{(iii)}. For $n\geq1$ and $x\in M_n(\mathscr{A})$,
\eqref{eq:extension-stinespring} gives
\[
\Phi^{(n)}(x)=\bigl((I_n\otimes N)\rho^{(n)}(x)\bigr)|_{\Hh^n},
\]
with $\Hh^n$ invariant for the product. If $x$ is subnormal, then
$\rho^{(n)}(x)$ is subnormal by \eqref{eq:BH}. It doubly commutes with
$I_n\otimes N$ because $\rho(\mathscr{A})$ is self-adjoint and commutes
with $N$. Lemma~\ref{lem:products} gives subnormality of the product and
hence of its invariant restriction. This proves \textup{(ii)}.

By contractivity of $\rho^{(n)}$ and minimality of $N$,
\[
\|\Phi^{(n)}(x)\|\leq\|N\|\,\|x\|=\|H\|\,\|x\|
\qquad(x\in M_n(\mathscr{A})).
\]
For $n=1$, this also gives boundedness under \textup{(iii)} when
$\mathscr{A}$ is unital. Since $\|H\|\leq\|\Phihat\|=\|\Phi\|$, the
estimate proves \eqref{eq:subnorm-stinespring}.
\end{proof}

\section*{Acknowledgments}
The authors thank the organizers of the 9th Linear Algebra Workshop
(Portoro\v{z}, Slovenia, June~2--6, 2025), during which the idea for
this paper originated.
The authors also thank OpenAI for providing access to ChatGPT through
the ChatGPT for Academic Researchers program.

\end{document}